 %\magnification=1200
\tolerance=10000
\hsize=15truecm \hoffset 0.4truecm           
\vsize=21truecm \voffset 1.5truecm

\def\Abstract#1\par{\centerline{\vbox{\hsize 11truecm
    \noindent {\bf Abstract.} #1\hfil}} \vskip 0.75cm}

\newcount\numsezione
\newcount\numcapitolo
\def\Capitolo #1{\advance \numcapitolo by 1 \numsezione=0
       \penalty -1000
        \bigskip
       {\bf \the \numcapitolo. \ #1 \hfill}
       \bigskip 
       \nobreak }

\def\Sezione #1{\ifnum\numsezione>0\goodbreak\fi
\noindent
       \global\advance \numsezione by 1
       \bigskip 
       \mark{ {\bf \the\numcapitolo.\the\numsezione}\enspace #1}
       {\bf \the \numcapitolo.\the \numsezione \ #1 \hfill}
       \nobreak\bigskip}

\def\References{\vskip1cm\goodbreak\noindent
    {\bf References \hfill}
    \vskip 0.25cm \nobreak}

\newcount \numproclaim

 \long\def\proclaim#1.#2 \par{ \advance\numproclaim by 1
         \goodbreak\ifdim\lastskip<\bigskipamount\removelastskip\bigskip\fi
         \noindent
        {\bf #1 \the \numproclaim.}
        {\sl #2\par}
         \bigskip}

    % MACRO PER ITEMS
\newdimen\myitemindent    %  indentazione, indipendente da
                                            %  \parindent, per l`enumerazione
                                             %  degli items

\myitemindent = 20 pt     %%%  o quanto preferisci
\def\myitem#1{\par\leavevmode\hbox to \myitemindent
               {\hbox{#1} \hfil}\ignorespaces}

\def\bull{{\vrule height.9ex width.8ex depth-.1ex}}

\def\k#1{{\cal #1}}

\def\lra{\longrightarrow}

\def\PP{\vbox{\hbox to 8.9pt{I\hskip-2.1pt P\hfill}}}
\def\II{\vbox{\hbox to 8.9pt{I\hskip-2.1pt I\hfill}}}

\font\tenmsbm=msbm10
\font\sevenmsbm=msbm7
\font\fivemsbm=msbm5
\newfam\amsBfam
\textfont\amsBfam=\tenmsbm
\scriptfont\amsBfam=\sevenmsbm
\scriptscriptfont\amsBfam=\fivemsbm
\def\bbb{\fam\the\amsBfam\tenmsbm}

\font\tenmsam=msam10
\font\sevenmsam=msam7
\font\fivemsam=msam5
\newfam\amsAfam
\textfont\amsAfam=\tenmsam
\scriptfont\amsAfam=\sevenmsam
\scriptscriptfont\amsAfam=\fivemsam
\def\aaa{\fam\the\amsAfam\tenmsam}

\mathchardef\gul"3\the\amsAfam52
\mathchardef\lug"3\the\amsAfam51
\mathchardef\R"5\the\amsBfam52

		\catcode`\"=12
	\font\kropa=lcircle10 scaled 1700
	\def\ybl{\setbox0=\hbox{\kropa \char"70} \kern1.5pt \raise.35pt \box0}
		\catcode`\"=\active

\def\bullet{}
	\def\normalsize{}

\catcode`\"=12
	
	\font\syx=cmsy10 scaled 1650

	\def\POUT{{\setbox0=\hbox{\syx\char"0E}\box0}}

\catcode`\"=\active

	\def\lra{\longrightarrow }

	\def\prod#1{\mathop{\POUT}\limits_{{\,#1\,}}}

\input pictex
\input dcpic.sty

%%%%%%%%%%%%%%%%%%%%%%%%%%%%%%%%%%%%%%%%%%%%%%%%%%
\numproclaim=0
\parindent=0pt

\centerline {\bf{Local algebra bundles and $\alpha$-jets of mappings.}}
\vskip 2\baselineskip 
\centerline {MARGHERITA BARILE$^*$}
\centerline{FIORELLA BARONE$^*$}
\centerline{WLODZIMIERZ M. TULCZYJEW$^{**}$}
\centerline{\hphantom{X}}
\centerline {$^*$Dipartimento di Matematica, Universit\`a di Bari, Italy.}
\centerline {$^{**}$Dipartimento di Fisica, Universit\`a di Camerino, Italy.}

\vskip 3\baselineskip
\Abstract We present the functor associated with a local algebra bundle and the differential structure of the double fibre bundle it produces when  applied to a differential manifold.

\vskip 2\baselineskip

\vskip 1\baselineskip
\noindent A.M.S. CLASSIFICATION: 13H10, 58A20

\Capitolo {Introduction and preliminaries.}

Weil algebras and related functors [14] have been quite extensively studied by a number of authors over the past decades ([2], [3], [4], [5], [7], [8], [12]). 
In the present note we generalize the method of associating a functor with a single local algebra. We start with a local algebra bundle $\alpha : A \to M$ on a differential manifold $M$ and define a functor which associates with a differential manifold $P$ a double fibre bundle
$$\vcenter{
\begindc{0}[3]
\obj(10,10){$M$}
\obj(30,18){$P$}
\obj(20,30){$AP$}
\mor(20,30)(10,10){$\alpha^s_P$}[\atright,\solidarrow]
\mor(20,30)(30,17){$\alpha^t_P$}
\enddc}\eqno(1)$$
The elements of $AP$ are called $\alpha$-jets on $P$.

In the first part of the paper, using some of the main results proved in [1] and [2], we give the proof of the differential structure of $(1)$.

In the second part we present a straightforward proof of the connection between the jets of mappings of Ehresmann type ([6], [7], [9], [10],  [11], [13]) and  a special class of $\alpha$-jets. 

\vskip1\baselineskip
A list of few basic definitions and conventions now follows.

Differential manifolds  will be assumed to be modelled on affine spaces and we will refer to the differential calculus on affine spaces presented in [1].

Let  $M$, $P$, $Q$ and $R$  be differential manifolds. 

\proclaim Definition.
A  {\it double bundle} with {\it total space} $Q$, {\it bases} $M$ and $P$ and  {\it standard fibre} $R$  is a diagram
$$\vcenter
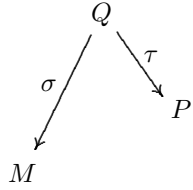
{
\begindc{0}[3]
\obj(10,10){$M$}
\obj(30,18){$P$}
\obj(20,30){$Q$}
\mor(20,30)(10,10){$\sigma$}[\atright,\solidarrow]
\mor(20,30)(30,17){$\tau$}
\enddc}$$
such that
\myitem {(i)} $\sigma : Q \to M$ and $\tau: Q \to P$ are surjective submersions;
\myitem {(ii)} for each $x\in M$ and $p\in P$, there are a chart $(\k U,\xi)$ around $x$, a chart $(\k W,\eta)$ around $p$ and a diffeomorphism
$$\chi : \sigma^{-1}(\k U)\cap \tau^{-1}(\k W) \lra \xi(\k U)\times \eta(\k W) \times R$$
such that the following  {\it  double local trivialization} diagram
$$\vcenter
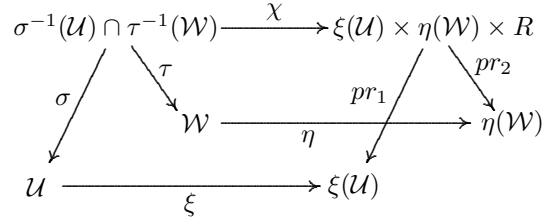
{
\begindc{0}[3]
\obj(10,10){$\k U$}
\obj(50,10){$\xi(\k U)$}
\obj(30,18){$\k W$}
\obj(70,18){$\eta(\k W)$}
\obj(20,30){$\sigma^{-1}(\k U)\cap \tau^{-1}(\k W)$}
\obj(60,30){$\xi(\k U)\times \eta(\k W) \times R$}
\mor(30,30)(50,30){$\chi$}
\mor(20,30)(10,10){$\sigma$}[\atright,\solidarrow]
\mor(60,30)(50,10){$pr_1$}[\atright,\solidarrow]
\mor(20,30)(30,17){$\tau $}
\mor(60,30)(70,17){$pr_2$}
\mor(10,10)(49,10){$\xi$}[\atright,\solidarrow]
\mor(30,18)(68,18){$\eta\phantom{xxxxx}$}[\atright,\solidarrow]
\enddc}$$
is commutative. 

A double bundle of the form
$$\vcenter
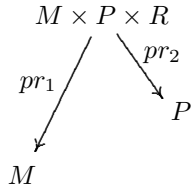
{
\begindc{0}[3]
\obj(10,10){$M$}
\obj(30,18){$P$}
\obj(20,30){$M\times P\times R$}
\mor(20,30)(10,10){$pr_1$}[\atright,\solidarrow]
\mor(20,30)(30,17){$pr_2$}
\enddc}$$
is called a {\it trivial double bundle}.

\proclaim Definition.
A {\it local algebra bundle} (l.a.b.)  is a (locally trivial) fibre bundle
$$\vcenter{
\begindc{0}[3]
\obj(10,10){$M$}
\obj(10,23){$A$}
\mor(10,23)(10,10){$\alpha$}
\enddc}$$
with a local algebra as its standard fibre.

For each $x\in M$, $A_x = \alpha^{-1}(x)$ is called the {\it local algebra over} $x$. 
If $I_{A_x}$ denotes the only maximal ideal of  $A_x$, we identify  $\R$ with $\R 1_{A_x}$, so that  $A_x = \R + I_{A_x}$.
Furthermore we denote by $0_{A_x}$ the only epimorphism from $A_x$ onto $\R$.

A  {\it local algebra bundle morphism} is a bundle morphism which induces a local algebra morphism on each fibre. 

\Capitolo {The functor associated with a local algebra bundle.}

Let
$$\vcenter{
\begindc{0}[3]
\obj(10,10){$M$}
\obj(10,23){$A$}
\mor(10,23)(10,10){$\alpha$}
\enddc}\eqno(2)$$
be a local algebra bundle.

Let $P$ be a differential manifold, and consider the set $A_{(x,\cdot)}P$ of all $\R$-algebra morphisms from $C^{\infty}(P)$ to $A_x $ which preserve the identity.
We set
$$AP=\bigcup_{x\in M}A_{(x,\cdot)}P.$$
The elements of $AP$ are called $\alpha$-{\it jets} on $P$.

We denote by 
$$\alpha^s_P : AP\lra M$$
the natural surjective mapping ({\it source}) such that for all $x\in M$,
$${(\alpha^s_P)}^{-1}(x) = A_{(x,\cdot)}P.$$
Let $u\in AP$ and $x = \alpha^s_P(u)$.
There is just one point $p\in P$ such that the composition $0_{A_x}\circ u$ is the evaluation mapping $ev_p$ at $p$ (cf., e.g., [7], p.296). 
So we have another surjective mapping ({\it target}) 
$$\alpha^t_P : AP\lra P$$
such that for all $p\in P$,
$$ (\alpha^t_P)^{-1}(p) = A_{(\cdot,p)}P =  \bigcup_{x\in M}\{ u\in A_{(x,\cdot)}P \, : \,  0_{A_x}\circ u = ev_p \}.\eqno(3)$$
We set
$$A_{(x,p)}P = A_{(x,\cdot)}P \cap A_{(\cdot,p)}P. $$
Then
$$A_{(x,\cdot)}P = \bigcup_{p\in P}A_{(x,p)}P \quad,\quad  A_{(\cdot,p)}P = \bigcup_{x\in M}A_{(x,p)}P$$

It is easy to prove the following

\proclaim Proposition. 
\endgraf
\myitem {(i)}  To a  l.a.b. $\alpha:A\to M$ corresponds a covariant functor
$$\vcenter{
\begindc{0}[3]
\obj(10,10){$M$}
\obj(30,18){$I$}
\obj(20,30){$A$}
\mor(20,30)(10,10){$\alpha^s$}[\atright,\solidarrow]
\mor(20,30)(30,17){$\alpha^t$}
\enddc}\eqno(4)$$
which associates with any differential manifold $P$ the diagram
$$\vcenter{
\begindc{0}[3]
\obj(10,10){$M$}
\obj(30,18){$P$}
\obj(20,30){$AP$}
\mor(20,30)(10,10){$\alpha^s_P$}[\atright,\solidarrow]
\mor(20,30)(30,17){$\alpha^t_P$}
\enddc}\eqno(5)$$
and with any differentiable mapping $\varphi: P\to Q$ the diagram
$$\vcenter{
\begindc{0}[3]
\obj(10,10){$M$}
\obj(50,10){$M$}
\obj(30,18){$P$}
\obj(70,18){$Q$}
\obj(20,30){$AP$}
\obj(60,30){$AQ$}
\mor(20,30)(60,30){$A\varphi$}
\mor(20,30)(10,10){$\alpha^s_P$}[\atright,\solidarrow]
\mor(60,30)(50,10){$\alpha^s_Q$}[\atright,\solidarrow]
\mor(20,30)(30,17){$\alpha^t_P$}
\mor(60,30)(70,17){$\alpha^t_Q$}
\mor(10,10)(50,10){$$}[\atright,\solidline]
\mor(10,9)(50,9){$$}[\atright,\solidline]
\mor(30,18)(70,18){$\varphi\phantom{xxxxx}$}[\atright,\solidarrow]
\enddc}$$
where, for all $u\in AP$,
$$A\varphi(u)= u \circ\varphi^*,\eqno(6)$$
i.e., for all $f\in C^{\infty}(Q)$,
$$A\varphi(u)(f)= u (f\circ\varphi).$$
\myitem{(ii)}  To a l.a.b. morphism
$$\vcenter{
\begindc{0}[3]
\obj(1,10){$M$}
\obj(30,10){$N$}
\obj(1,25){$A$}
\obj(30,25){$B$}
\mor(1,25)(30,25){$\kappa$}
\mor(1,10)(30,10){$\mu$}
\mor(1,25)(1,10){$\alpha$}[\atright,\solidarrow]
\mor(30,25)(30,10){$\beta$}
\enddc}$$
corresponds the functor
$$\vcenter{
\begindc{0}[1]
\obj(39,13){$N$}
\obj(143,39){$I$}
\obj(91,65){$B$}
\obj(13,78){$M$}
\obj(117,104){$I$}
\obj(65,130){$A$}
\mor(65,130)(13,78){$\alpha^s$}[\atright,\solidarrow]
\mor(65,130)(117,104){$\alpha^t$}
\mor(91,65)(39,13){$\beta^s$}[\atright,\solidarrow]
\mor(91,65)(143,39){$\beta^t$}
\mor(13,78)(39,13){$\mu$}[\atright,\solidarrow]
\mor(65,130)(91,65){$\kappa$}
\mor(117,104)(143,39){$$}[\atright,\solidline]
\mor(120,104)(146,39){$$}[\atright,\solidline]
\enddc}$$
which associates with a differential manifold $P$ the diagram
$$\vcenter{
\begindc{0}[1]
\obj(39,13){$N$}
\obj(143,39){$P$}
\obj(91,65){$BP$}
\obj(13,78){$M$}
\obj(117,104){$P$}
\obj(65,130){$AP$}
\mor(65,130)(13,78){$\alpha^s_P$}[\atright,\solidarrow]
\mor(65,130)(117,104){$\alpha^t_P$}
\mor(91,65)(39,13){$\beta^s_P$}[\atright,\solidarrow]
\mor(91,65)(143,39){$\beta^t_P$}
\mor(13,78)(39,13){$\mu$}[\atright,\solidarrow]
\mor(65,130)(91,65){$\kappa_P$}
\mor(117,104)(143,39){$$}[\atright,\solidline]
\mor(120,104)(146,39){$$}[\atright,\solidline]
\enddc}$$
where
$$\kappa_P(u)=\kappa\circ u,$$
and with a differentiable mapping $\varphi: P\to Q$ the diagram
$$\vcenter
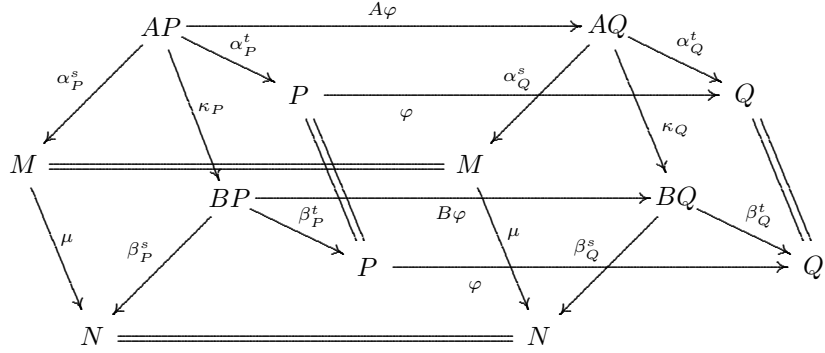
{
\begindc{0}[1]
\obj(39,13){$N$}
\obj(143,39){$P$}
\obj(91,65){$BP$}
\obj(13,78){$M$}
\obj(117,104){$P$}
\obj(65,130){$AP$}
\obj(208,13){$N$}
\obj(312,39){$Q$}
\obj(260,65){$BQ$}
\obj(182,78){$M$}
\obj(286,104){$Q$}
\obj(234,130){$AQ$}
\obj(260,91){${}_{\kappa_Q}$}
\mor(234,130)(182,78){${}_{\alpha^s_Q}$}[\atright,\solidarrow]
\mor(234,130)(286,104){${}_{\alpha^t_Q}$}
\mor(260,65)(208,13){${}_{\beta^s_Q}$}[\atright,\solidarrow]
\mor(260,65)(312,39){${}_{\beta^t_Q}$}
\mor(182,78)(208,13){${}^{\mu}$}
\mor(234,130)(260,65){$$}
\mor(286,104)(312,39){$$}[\atright,\solidline]
\mor(289,104)(315,39){$$}[\atright,\solidline]
\mor(65,130)(13,78){${}_{\alpha^s_P}$}[\atright,\solidarrow]
\mor(65,130)(117,104){${}_{\alpha^t_P}$}
\mor(91,65)(39,13){${}_{\beta^s_P}$}[\atright,\solidarrow]
\mor(91,65)(143,39){${}_{\beta^t_P}$}
\mor(13,78)(39,13){${}_{\mu}$}
\mor(65,130)(91,62){${}_{\kappa_P}$}
\mor(117,104)(143,39){$$}[\atright,\solidline]
\mor(120,104)(146,39){$$}[\atright,\solidline]
\mor(39,13)(208,13){$$}[\atright,\solidline]
\mor(39,11)(208,11){$$}[\atright,\solidline]
\mor(143,39)(312,39){${}_{\varphi}\phantom{xxxxxxxxxxxxxxx}$}[\atright,\solidarrow]
\mor(91,65)(260,65){${}_{B\varphi}$}[\atright,\solidarrow]
\mor(13,78)(182,78){$$}[\atright,\solidline]
\mor(13,76)(182,76){$$}[\atright,\solidline]
\mor(117,104)(286,104){${}_{\varphi}\phantom{xxxxxxxxxxxxxxx}$}[\atright,\solidarrow]
\mor(65,130)(234,130){${}_{A\varphi}$}
\enddc}$$
\hfill$\bull$

\vskip0.5\baselineskip
In this section we will prove that diagram $(5)$ is a double bundle. 

To this end, we first consider the functor associated with a trivial l.a.b. over an affine space and apply it to another affine space.
The resulting trivial double bundle will be the model for the  structure we are looking for.

In the sequel, $E$ will denote an affine space modelled on a vector space $X$ and $G$ an affine space modelled on a vector space $Z$.

\proclaim Lemma.
Let
$$\vcenter{
\begindc{0}[3]
\obj(10,10){$E$}
\obj(10,23){$E\times \bar A$}
\mor(10,23)(10,10){$\varepsilon=pr_1$}
\enddc}\eqno(7)$$
be a trivial l.a.b. over the  affine space $E$.
The corresponding functor, when applied to the affine space $G$, yields a trivial double bundle with an affine space as its total space.

{\it Proof}.
By definition, when we apply the functor associated with l.a.b. (7) to $G$, we get
$$\vcenter{
\begindc{0}[3]
\obj(10,10){$E$}
\obj(30,18){$G$}
\obj(20,30){$(E\times \bar A)G$}
\mor(20,30)(10,10){$\varepsilon^s_G$}[\atright,\solidarrow]
\mor(20,30)(30,17){$\varepsilon^t_G$}
\enddc}\eqno(8)$$
Then we have 
$$\eqalign{
(E\times \bar A)G &= \bigcup_{x\in E}(\varepsilon^s_G)^{-1}\{x\}\cr
                                &= \bigcup_{x\in E}\big(\{x\}\times \bar A\big)G\cr
                                &= E\times \bar A G. \cr
}\eqno(9)$$
We have proved the following identification (cf. [2], Prop.9)
 $$\bar AG = G\times(I_{\bar A}\otimes Z), \eqno(10)$$ 
 where $I_{\bar A}$ is the maximal ideal of $\bar A$.
As a consequence we have
$$(E\times \bar A)G= E\times G\times(I_{\bar A}\otimes Z),\eqno(11)$$
so that $(E\times \bar A)G$ is an affine space and diagram (8) becomes the following trivial double bundle
$$\vcenter{
\begindc{0}[3]
\obj(10,10){$E$}
\obj(30,18){$G$}
\obj(20,30){$ E\times G\times (I_{\bar A}\otimes Z)$}
\mor(20,30)(10,10){$pr_1$}[\atright,\solidarrow]
\mor(20,30)(30,17){$pr_2$}
\enddc}$$
\hfill$\bull$

\vskip0.5\baselineskip
Let us notice that  equalities $(10)$ and $(11)$ still hold if $E$ and $G$ are replaced by open subsets owing to  the local character of $\alpha$-jets which we briefly recall.

Let $\k W$ be an open subset of $P$.
From the action of the functor (4) on the inclusion mapping
$$\iota_{\k W}:\k W\to P$$
we obtain the mapping
$$A\iota_{\k W} : A\k W\lra AP$$
whose restriction to any fibre $A_{(\cdot,p)}\k W$ of  $\alpha_P^t$
$$A_p\iota_{\k W} : A_{(\cdot,p)}\k W\lra A_{(\cdot,p)}P$$
is a bijection (cf. [2], Sec.4.2).
As a consequence, owing to the identification
$$ A_{(\cdot,p)}\k W =  A_{(\cdot,p)}P,$$
we have that $A\k W$ is a subset of $AP$ (local character of $\alpha$-jets) and that $A\iota_{\k W}$ is the inclusion mapping $\iota_{A\k W}$.
This is illustrated by the following diagram
\endgraf
\vskip0.5\baselineskip
$$\vcenter{
\begindc{0}[3]
\obj(10,10){$M$}
\obj(40,10){$M$}
\obj(35,18){$\k W$}
\obj(65,18){$P$}
\obj(20,30){$A\k W$}
\obj(50,30){$AP$}
\mor(20,30)(50,30){$\iota_{A\k W}$}[\atleft,\solidarrow]
\mor(20,30)(10,10){$\alpha^s_{\k W}$}[\atright,\solidarrow]
\mor(50,30)(40,10){$\alpha^s_P$}[\atright,\solidarrow]
\mor(20,30)(35,18){$\alpha^t_{\k W}$}
\mor(50,30)(65,18){$\alpha^t_P$}
\mor(10,10)(40,10){$$}[\atright,\solidline]
\mor(10,11)(40,11){$$}[\atright,\solidline]
\mor(35,18)(65,18){$\iota_{\k W}$}[\atright,\solidarrow]
\enddc}$$

\vskip0.5\baselineskip
The following step is to apply the functor associated with the l.a.b. morphism given by a trivialization of $(2)$ to a chart of differential manifold $P$. The result  will be a chart of $AP$.

In the sequel, $M$ will denote a differential manifold modelled on the affine space $E$ and $P$ a differential manifold modelled on the affine space $G$.

Let
$$\vcenter{
\begindc{0}[3]
\obj(10,10){$\k U$}
\obj(45,10){$ \xi(\k U)$}
\obj(10,25){$A_{\xi}=\alpha^{-1}(\k U)$}
\obj(45,25){$ \xi(\k U)\times  \bar A$}
\mor(18,25)(41,25){$\xi^A$}
\mor(10,10)(43,10){$\xi$}[\atright,\solidarrow]
\mor(10,25)(10,10){$\alpha$}[\atright,\solidarrow]
\mor(45,25)(45,10){$\varepsilon=pr_1$}
\enddc}\eqno(12)$$
be a local trivialization  of l.a.b. (2) with  $\xi : \k U \to \xi(\k U)$ an admissible chart on $M$ and $\bar A$  the standard fibre.

If we apply the functor corresponding to  l.a.b. isomorphism $(12)$ to an admissible chart $(\eta ,\k W)$ on $P$, we obtain  the diagram

\vskip0.5\baselineskip
$$\vcenter{
\begindc{0}[1]
\obj(39,13){$H$}
\obj(143,39){$\k W$}
\obj(91,65){$( H\times\bar A)\k W$}
\obj(13,78){$\k U$}
\obj(117,104){$\k W$}
\obj(65,130){$A_{\xi}\k W$}
\obj(208,13){$H$}
\obj(312,39){$K$}
\obj(260,65){$(H\times\bar A)K$}
\obj(182,78){$\k U$}
\obj(286,104){$K$}
\obj(234,130){$A_{\xi}K$}
\obj(260,91){${\xi^A_K}$}
\mor(234,130)(182,78){${\alpha^s_K}$}[\atright,\solidarrow]
\mor(234,130)(286,104){${\alpha^t_K}$}
\mor(260,65)(208,13){${\varepsilon^s_K}$}[\atright,\solidarrow]
\mor(260,65)(312,39){${\varepsilon^t_K}$}
\mor(182,78)(208,13){${\xi}$}
\mor(234,130)(260,63){$$}
\mor(286,104)(312,39){$$}[\atright,\solidline]
\mor(289,104)(315,39){$$}[\atright,\solidline]
\mor(65,130)(13,78){${\alpha^s_{\k W}}$}[\atright,\solidarrow]
\mor(65,130)(117,104){${\alpha^t_{\k W}}$}
\mor(91,65)(39,13){${\varepsilon^s_{\k W}}$}[\atright,\solidarrow]
\mor(91,65)(143,39){${\varepsilon^t_{\k W}}$}[\atright,\solidarrow]
\mor(13,78)(39,13){${\xi}$}
\mor(65,130)(91,63){${\xi^A_{\k W}}$}
\mor(117,104)(143,39){$$}[\atright,\solidline]
\mor(120,104)(146,39){$$}[\atright,\solidline]
\mor(39,13)(208,13){$$}[\atright,\solidline]
\mor(39,11)(208,11){$$}[\atright,\solidline]
\mor(143,39)(312,39){${\eta}\phantom{xxxxxxxxxxxxxxx}$}[\atright,\solidarrow]
\mor(105,65)(246,65){${(H\times\bar A)\eta\phantom{xxxx}}$}[\atright,\solidarrow]
\mor(13,78)(182,78){$$}[\atright,\solidline]
\mor(13,76)(182,76){$$}[\atright,\solidline]
\mor(117,104)(286,104){${\eta}\phantom{xxxxxxxxxxxxxxx}$}[\atright,\solidarrow]
\mor(70,132)(234,132){${A_{\xi}\eta}$}
\enddc}\eqno(13)$$

where, for the sake of simplicity, we have set $H = \xi(\k U)$ and $K = \eta(\k W)$.

From diagram (13) we can select the subdiagram

$$\vcenter{
\begindc{0}[2]
\obj(10,10){$\k U$}
\obj(50,10){$ \xi(\k U)$}
\obj(30,20){$\k W$}
\obj(70,20){$ \eta(\k W)$}
\obj(20,40){$A_{\xi}\k W$}
\obj(70,40){$\big(\xi(\k U)\times {\bar  A}\big) \eta(\k W)$}
\mor(20,40)(10,10){$\alpha^s_{\k W}$}[\atright,\solidarrow]
\mor(20,40)(30,18){$\alpha^t_{\k W}$}
\mor(60,40)(50,10){$ \varepsilon^s_{ \eta(\k W)}$}[\atright,\solidarrow]
\mor(60,40)(70,18){$ \varepsilon^t_{ \eta(\k W)}$}
\mor(10,10)(49,10){$\xi$}[\atright,\solidarrow]
\mor(30,20)(68,20){$\eta\phantom{xxx}$}[\atright,\solidarrow]
\mor(23,40)(55,40){$A_{(\xi,\eta)}$}
\enddc}\eqno(14)$$
where
$$A_{(\xi,\eta)}:=\xi_{ \eta(\k W)}^A\circ A_{\xi}\eta$$
is a bijection, which, owing to Lemma 4, is a chart of $AP$.

\proclaim Proposition.
$AP$ is a differential manifold.

{\it Proof.}
$AP$ will be a differential manifold modelled on the affine space $(E\times\bar A) G$ provided that the charts $A_{(\xi,\eta)}$ are $C^{\infty}$-compatible.

Without loss of generality we will consider charts $\xi'$ and $\eta'$ on $M$ and $P$ having the same domains $\k U$ and $\k W$ as $\xi$ and $\eta$, respectively.
It follows that $A_{\xi} = A_{\xi'}$.

We have to prove that the composition
$$\eqalign{
A_{(\xi',\eta')} &\circ A_{(\xi,\eta)}^{-1} =
{\xi'}_{ \eta'(\k W)}^A\circ A_{\xi'}\eta' \circ  (A_{\xi}\eta)^{-1} \circ (\xi_{ \eta(\k W)}^A)^{-1}\cr
&= {\xi'}_{ \eta'(\k W)}^A\circ A_{\xi}(\eta' \circ  \eta^{-1}) \circ (\xi_{ \eta(\k W)}^A)^{-1} :  \big(\xi(\k U)\times  \bar A\big)\eta(\k W)\lra \big( \xi'(\k U)\times  \bar A\big)\eta'(\k W)\cr 
}\eqno(15)$$
is  differentiable.

\myitem{(i)} Let us consider the transition function obtained from $(12)$ and the analogous local trivialization of (2) corresponding to $\xi'$:

$$\vcenter{
\begindc{0}[3]
\obj(10,10){$\xi(\k U)$}
\obj(45,10){$ \xi'(\k U)$}
\obj(10,25){$ \xi(\k U)\times  \bar A$}
\obj(45,25){$ \xi'(\k U)\times  \bar A$}
\mor(18,25)(41,25){$\xi'^A \circ (\xi^A)^{-1}$}
\mor(10,10)(43,10){$\xi'\circ \xi^{-1}$}[\atright,\solidarrow]
\mor(10,25)(10,10){$pr_1$}[\atright,\solidarrow]
\mor(45,25)(45,10){$pr_1$}
\enddc}\eqno(16)$$

It is a differentiable l.a.b. morphism.
It follows that for all $\xi(x) \in \xi(\k U)$ the induced mapping 
$$
{\xi'}^A\big\vert_{A_x } \circ    (\xi^A\big\vert_{A_x })^{-1} : \xi(x)\times\bar A \lra \xi'(x)\times\bar A\; : \;
\big(\xi(x),a\big) \mapsto \big(\xi'\circ\xi^{-1}(\xi(x)), \Xi^{\xi(x)}(a)\big) $$
determines a local algebra automorphism
$$\Xi^{\xi(x)} : \bar A \lra \bar A\eqno(17)$$
and, moreover, that the mapping $ \xi(x) \mapsto \Xi^{\xi(x)}$ is  differentiable.

Owing to Proposition 7 in [2], when we apply the functor corresponding to (17) to $\eta'(\k W)$ we get the following map:
$$\Xi^{\xi(x)} _{\eta'(\k W)} : \bar A \eta'(\k W) \lra  \bar A \eta'(\k W) \; : \; \bar u \mapsto\Xi^{\xi(x)}  \circ \bar u.\eqno(18)$$

The above mapping is differentiable since, by Proposition 9  in [2], for any choice of a Cartesian coordinate system $(y^1,\dots,y^m)$ in $G$, it admits the following global coordinate expression

$$\eta'(\k W)\times({\underbrace{I_{\bar A}\times\dots\times I_{\bar A}}_{m\ {\rm \scriptstyle times}}})  
\lra  \eta'(\k W)\times({\underbrace{I_{\bar A}\times\dots\times I_{\bar A}}_{m\ {\rm \scriptstyle times}}})$$ 
$$\quad\quad\quad \big(\eta'(p), (\bar u^j)_{j=1\dots,m}) \mapsto  \big(\eta'(p),(\Xi^{\xi(x)} (\bar u^j))_{j=1\dots,m}\big), \eqno(19)$$
where we have set $\bar u^j = \bar u(y^j) -  \bar u(\eta'(p))$ for all $j=1\dots,m$.

\myitem{(ii)} Now we can prove that the mapping (15) is differentiable. It acts as follows:
$$\vcenter{
\begindc{0}[2]
\obj(105,20){$ {\xi'}^A\big\vert_{A_x } \circ    (\xi^A\big\vert_{A_x })^{-1} \circ v \circ   (\eta'\circ\eta^{-1})^*$}
\obj(94,30){$ (\xi^A\big\vert_{A_x })^{-1} \circ v \circ  (\eta'\circ\eta^{-1})^*$}
\obj(80,40){$ (\xi^A\big\vert_{A_x })^{-1} \circ v $}
\obj(1,40){$ v $}
\mor(10,40)(65,40){${(\xi_{ \eta(\k W)}^A)^{-1} }$}
\mor(10,30)(65,30){${A_{\xi}(\eta'\circ \eta^{-1}) }$}
\mor(10,20)(65,20){${{\xi'}_{ \eta'(\k W)}^A }$}
\enddc}\eqno(20)$$

Owing to $(9)$, we can use the following identification
$$\eqalign{
\big( \xi(\k U)\times  \bar A\big)\eta(\k W) &\cong  \xi(\k U)\times  \bar A \eta(\k W)\cr
v & \cong \big( \xi(x),\bar v\big)\cr
}$$
where $\bar v = pr_2 \circ v$.
With the above identification, $(20)$ becomes

$$\eqalign{
\xi(\k U)\times  \bar A \eta(\k W) \lra & \xi'(\k U)\times  \bar A \eta'(\k W)\cr
\big( \xi(x),\bar v\big) \lra & \big(\xi'(x) , \Xi^{\xi(x)}\circ \bar v \circ (\eta'\circ\eta^{-1})^*\big)\cr
& = \Big(\xi'\circ\xi^{-1}\big( \xi(x)\big), \Xi^{\xi(x)}_{\eta'(\k W)}[\bar A (\eta'\circ\eta^{-1})(\bar v)]\Big)\cr
}$$
Owing to Proposition 10 in [2] and (i), we conclude that  (15) is a differentiable mapping.
\hfill $\bull$

\vskip1\baselineskip
Up to now we have proved that $AP$ is a differential manifold.

The next step is to show that both $\alpha_P^s$ and $\alpha_P^t$ are surjective submersions.

\proclaim Proposition.
The diagram
$$\vcenter{
\begindc{0}[3]
\obj(10,10){$M$}
\obj(10,23){$AP$}
\mor(10,23)(10,10){$\alpha^s_P$}
\enddc}$$
is a surjective submersion.

{\it Proof.}
From diagram $(14)$ we select the subdiagram
$$\vcenter{
\begindc{0}[2]
\obj(10,10){$\k U$}
\obj(50,10){$\xi(\k U)$}
\obj(20,40){$A_{\xi}\k W$}
\obj(90,40){$\big(\xi(\k U)\times\bar A\big) \eta(\k W) = \xi(\k U)\times \bar A\eta(\k W)$}
\mor(20,40)(10,10){$\alpha^s_{\k W}$}[\atright,\solidarrow]
\mor(60,40)(50,10){$ \varepsilon^s_{\eta(\k W)}= pr_1$}%[\atright,\solidarrow]
\mor(10,10)(49,10){$\xi$}
\mor(22,40)(56,40){$A_{(\xi,\eta)}$}
\enddc}$$

It is the coordinate expression of $\alpha^s_P$, which is a projection, and hence a surjective submersion.

\hfill$\bull$

\proclaim Proposition.
The diagram
$$\vcenter{
\begindc{0}[3]
\obj(10,10){$P$}
\obj(10,23){$AP$}
\mor(10,23)(10,10){$\alpha^t_P$}
\enddc}$$
is a surjective submersion.

{\it Proof.}
From diagram $(14)$ we select the subdiagram
$$\vcenter{
\begindc{0}[2]
\obj(30,20){$\k W$}
\obj(70,20){$\eta(\k W)$}
\obj(20,40){$A_{\xi}\k W$}
\obj(90,40){$\big(\xi(\k U)\times\bar A\big) \eta(\k W) = \xi(\k U)\times \bar A\eta(\k W)$}
\mor(20,40)(30,20){$\alpha^t_{\k W}$}
\mor(60,40)(70,20){$ \varepsilon^t_{\eta(\k W)}=\bar\alpha^t_P\circ  pr_2$}
\mor(30,20)(68,20){$\eta\phantom{xxx}$}
\mor(22,40)(58,40){$A_{(\xi,\eta)}$}
\enddc}$$

It is the coordinate expression of $\alpha^t_P$, which is a projection, and hence a surjective submersion.
\hfill$\bull$

\proclaim Theorem.
Diagram 
$$\vcenter{
\begindc{0}[3]
\obj(10,10){$M$}
\obj(30,18){$P$}
\obj(20,30){$AP$}
\mor(20,30)(10,10){$\alpha^s_P$}[\atright,\solidarrow]
\mor(20,30)(30,17){$\alpha^t_P$}
\enddc}$$
 is a double bundle.

{\it Proof.}
We only have to prove part (ii) of Definition 1.

Let $(\k U,\xi)$ be a chart around a point $x$ of $M$ and $(\k W,\eta)$ a chart around a point $p$ of $P$.
We have 

$$\eqalign{
(\alpha^s_P)^{-1}(\k U) \bigcap (\alpha^t_P)^{-1}(\k W) & =
\left(\bigcup_{x\in \k U}A_{(x,\cdot)}P\right) \bigcap \left( \bigcup_{p\in \k W}A_{(\cdot,p)}P\right)\cr
& =\left( \bigcup_{x\in \k U}\bigcup_{p\in P}A_{(x,p)}P\right) \bigcap \left( \bigcup_{p\in \k W} \bigcup_{x\in M}A_{(x,p)}P\right)\cr
& =\bigcup_{x\in \k U} \bigcup_{p\in \k W}A_{(x,p)}P  = \bigcup_{x\in \k U} \bigcup_{p\in \k W}A_{(x,p)}\k W \cr
& = \bigcup_{x\in \k U} A_{(x,\cdot)}\k W = \big(\alpha^{-1}(\k U)\big)\k W\cr
&= A_{\xi}\k W.
}$$

Therefore diagram $(14)$ is a double local trivialization of $(\alpha^s_P,\alpha^t_P)$.
\hfill$\bull$

\Sezione {$\alpha$-jets of mappings.}

In this section we present  the connection between the jets of mappings of Ehresmann type  and  a special class of $\alpha$-jets. 

Let $x\in M$.
In the algebra $C^{\infty}(M)$ consider the maximal ideal $I_0(M,x)$ of all functions vanishing at $x$, and
its powers
$$I^k(M,x) = \big( I_0(M,x) \big)^{k+1}\quad\quad \forall k\in\bbb{N}.$$

We can use these ideals in two ways.

\myitem {(i)} Consider the space $C^{\infty}(M|P)$ of all differentiable mappings from $M$ to  $P$ and
introduce, in the Cartesian product $C^{\infty}(M|P) \times M$, the following equivalence relation:
$$(\varphi',x')\sim_k (\varphi, x) \iff
    \cases{x' = x \cr
    g \circ\varphi' - g\circ\varphi  \in I^k(M,x) \quad\quad \forall g \in  C^{\infty}(P)    
   }.$$
The equivalence class of $(\varphi, x)$ is denoted by $j^k\varphi(x)$  and called the $k$-jet of
$\varphi$ at $x$. The set of all $k$-jets is denoted by  $J^k(M|P)$.

\myitem {(ii)} Consider the special case of $P=\R$.
For each $x\in M$,  the set of $k$-jets of functions on $M$ at $x$, denoted by $A^k_x(M)$, can be identified (cf. [7], Remark 12.7) with the local algebra
$$A^k_x(M) =  \R + I_0(M,x)/I^k(M,x),$$
and then we have the local algebra bundle
$$\vcenter{
\begindc{0}[3]
\obj(10,10){$M$}
\obj(10,30){$A^k(M)$}
\mor(10,30)(10,10){$\alpha^k$}[\atright,\solidarrow]
\enddc}\eqno(21)$$

Now, if we apply the functor associated with (21) to differential manifold $P$, we obtain the double bundle

$$\vcenter
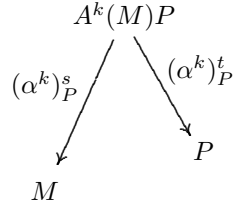
{
\begindc{0}[3]
\obj(10,8){$M$}
\obj(30,13){$P$}
\obj(20,30){$A^k(M)P$}
\mor(20,30)(10,8){$(\alpha^k)^s_P$}[\atright,\solidarrow]
\mor(20,30)(30,12){$(\alpha^k)^t_P$}
\enddc}$$

Constructions (i) and (ii) are, actually, closely related, as shown in the following theorem.

\proclaim Theorem.
The mapping
$$\chi : J^k(M|P) \lra A^k(M)P$$
 defined by
$$\chi( j^k\varphi(x)) : g \mapsto  j^k(g\circ\varphi)(x)$$
is a bijection.

Proof.
Note that $\chi$ is well defined  by definition of  $J^k(M|P)$.
We prove bijectivity by showing that for every $u\in A^k(M)P$ there is a unique $j^k\varphi(x)\in
J^k(M|P)$ such that
$$\chi\big(j^k\varphi(x)\big) = u.$$
Let $p = (\alpha^k)^t_P(u)$ and $\eta$ be any chart  on $P$ around $p$.
Given a Cartesian frame $\phi$ on the $n$-dimensional affine space $G$ centred at $\eta(p)$, put
$$\zeta = \phi\circ\eta$$
and
$$\zeta^h = x^h\circ\zeta,$$
where $x^h$ is the $h$-th natural coordinate function on $\R^n$.
For all  indices $h$ let $\varphi^h\in C^{\infty}(M)$ be such that
$$u(\zeta^h) = j^k\varphi^h(x).$$
Then, owing to (3),
$$\varphi^h(x) = \zeta^h(p).$$
Up to replacing  $\varphi^h$ with another representative of the same jet, we may assume that
$\varphi^h$ is not constant. 
Since $M$ is connected, it follows that the image of $\varphi^h$ is an
interval in $\R$. 
As a consequence the image of
$$\varphi_{\zeta} = (\varphi^1, \dots , \varphi^n)$$
contains an open neighborhood of
$$\zeta(p) = \varphi_{\zeta} (x).$$
The composition $\zeta^{-1}\circ  \varphi_{\zeta}$ defines a mapping $\varphi\in C^{\infty}(M|P)$ such that
$$ \varphi^h = \zeta^h \circ  \varphi.$$
We will prove that the mapping $\varphi$ is a representative of the required jet in $J^k(M|P)$.

Let $g\in  C^{\infty}(P)$, and let $g_{\zeta}=g\circ\zeta^{-1}$.
Owing to the local character of $u$, we have 
$$u(g) = u(g_{\zeta}\circ\zeta).$$
We will write, for $g_{\zeta}$, the Taylor formula at $\zeta(p)=0$, of order $k$ 
$$g_{\zeta}= {\k P}^k(x^1,\dots,x^n) + R^k$$
with Lagrange remainder
$$R^k = {\k P}^{k+1}H,$$
where ${\k P}^{k+1}$ is a homogeneous polynomial of degree $k+1$ at $0$.

Then we have
$$\eqalign{
g_{\zeta}\circ\zeta &=  {\k P}^k(x^1\circ\zeta,\dots,x^n\circ\zeta) + R^k\circ\zeta \cr
&=  {\k P}^k(\zeta^1,\dots, \zeta^n) +({\k P}^{k+1}H)\circ\zeta .
}$$
As a consequence
$$u(g) =  u\big({\k P}^k(\zeta^1,\dots, \zeta^n)\big) + u({\k P}^{k+1}\circ\zeta) u(H\circ\zeta),$$
where
$$u({\k P}^{k+1}\circ\zeta)=u({\k P}^{k+1}(\zeta^1, \dots, \zeta^n))=
{\k P}^{k+1}(u(\zeta^1), \dots,u( \zeta^n)) = 0 \eqno(22)$$
since $u( \zeta^h)\in I_0(M,x)/I^k(M,x)$. 
Therefore
$$\eqalign{
u(g) &= {\k P}^k\big(u(\zeta^1),\dots, u(\zeta^n)\big) \cr
&={\k P}^k\big(j^k\varphi^1(x),\dots, j^k\varphi^n(x)\big)  \cr
&= j^k\big( {\k P}^k\big(\varphi^1,\dots, \varphi^n\big)\big)(x)  \cr
&= j^k\big( {\k P}^k\big(\zeta^1\circ \varphi,\dots,\zeta^n\circ \varphi\big)\big)(x) \cr
&= j^k\big( {\k P}^k\big(\zeta^1,\dots,\zeta^n\big)\circ \varphi\big)(x)  \cr
&= j^k\big( {\k P}^k\big(\zeta^1,\dots,\zeta^n\big)\circ \varphi\big)(x).  \cr
}$$
Due to the local character of jets, we have
$$\eqalign{
u(g) &= j^k\big( {\k P}^k\big(\zeta^1,\dots,\zeta^n\big)\circ \varphi\big)(x) \cr
&= j^k\big(g_{\zeta}\circ\zeta \circ \varphi - R^k\circ\zeta\circ\varphi\big)(x)  \cr
&= j^k(g \circ \varphi)(x) - j^k(R^k\circ\zeta\circ\varphi\big)(x) \cr
&= j^k(g \circ \varphi)(x). \cr
}$$

The last equality is justified by arguments similar to those used in (22).

The uniqueness of $j^k\varphi(x)$ trivially follows by definition.
\hfill$\bull$

 \References
 
 [1] M. Barile, F. Barone and W.M. Tulczyjew, {\it Polarizations and Differential Calculus in Affine Spaces}, Linear Multilinear Algebra. {\bf 55}, 121--146 (2007).

 [2] M. Barile, F. Barone and W.M. Tulczyjew, {\it The category of local algebras and points proches}, Bull. Math. Soc. Sci. Math. R. S. Roumanie (N.S.),  {\bf 50}, 3--31 (2007). 
 
 [3] E. Dubuc,  {\it  Sur les mod\`eles de la g\'eom\'etrie diff\'erentielle synth\'etique}, Cah. Topol. G\'eom. Diff\'er. Cat\'eg., {\bf XX}, 234--278 (1979).
 
 [4] E. Dubuc and J. Zilber,  {\it  Weil Prolongations of Banach Manifolds in an Analytic Model of SDG}, Cah. Topol. G\'eom. Diff\'er. Cat\'eg., {\bf XLVI}, 83--98 (2005).

 [5] D.J. Eck, {\it  Product--Preserving Functors on Smooth Manifolds}, J. Pure Appl. Algebra, {\bf 42}, 133--140 (1986).

 [6] C. Ehresmann, {\it  Les prolongements d'une vari\'et\'e diff\'erentiable}, C.R. Acad. Sci., {\bf 233}, 598--600, 777--779 (1951).

 [7] I. Kol\'a\v{r}, P.W. Michor and J. Slov\'ak, {\it Natural Operations in Differential Geometry}, Springer (1993).
 
 [8] I.ÊKol\'ar, Ê{\it Bundle functors of the jet type}, in: I.ÊKol\'ar et al. (eds.), Differential geometry and applications.  Proceedings of the 7th international conference, DGA 98, and satellite conference of ICM in Berlin, Brno, Czech Republic, August 10-14, 1998. Brno: Masaryk University, 231--237 (1999).

[9] P. Liebermann, {\it Introduction to the theory of semi-holonomic jets}, Arch. Math. Brno, {\bf 33}, 173--189 (1996).

[10] G. Pidello, W.M. Tulczyjew, {\it Derivations of differential forms on jet bundles}, Ann. Mat. Pura Appl., IV. Ser., {\bf CXLVII}, 249--352 (1987).

[11] J.  Pradines, {\it  Fibr\'es vectoriels doubles et calcul des jets non holonomes}, Esquisses Math., {\bf 29}, Amiens (1977).

[12] V.V. Shurygin, {\it The structure of smooth mappings over Weil algebras and the category of manifolds over algebras}, Lobachevskii J. Math., {\bf 5}, 29-55 (1999).

[13] W.M. Tulczyjew, {\it Evolution of Ehresmann's jet theory}. To appear in: Banach Center Publications, Institute of Mathematics, Polish Academy of Sciences, Warszawa. 

[14] A. Weil, {\it Th\'eorie des points proches sur les vari\'et\'es diff\'erentiables}, Colloques internat. Centre nat. Rech. Sci., {\bf 52}, 111--117 (1953).

 \bye